\documentclass[10pt]{article}
\usepackage{amssymb}
\usepackage{amsthm}\usepackage{epsf,epsfig,amsfonts}
\usepackage{amsmath}
\usepackage{graphicx} 
\newcommand{\wt}[1]{\widehat{#1}}

\newcommand{\tg}{\hat g}
\usepackage{a4wide}

\newtheorem{corollary*}{Corollary}

\newcommand{\be}{\begin{equation}}
\newcommand{\ee}{\end{equation}}

\newcommand{\R}{\mathbb{R}}

\newcommand{\weg}[1]{}

\usepackage{epsf,epsfig,amsfonts,a4wide}

\usepackage{amsfonts}
\usepackage{amsmath,amssymb,amsthm}
\usepackage{url}
\usepackage{amsmath,amssymb,amsthm}

\usepackage{amsfonts}
\usepackage{amsmath,amssymb,amsthm}
\usepackage{url}

\newtheorem{Th}{Theorem}

\newtheorem{Lemma}{Lemma}
\newtheorem{Cor}{Corollary}

\theoremstyle{remark}

\newtheorem{Rem}{Remark}

\newcommand{\Id}{\mbox{\rm Id}}

\newcommand{\const}{\mbox{\rm const}}
\title{There exist no 4-dimensional geodesically equivalent metrics with the same stress-energy tensor}
\date{} \author{ Volodymir  Kiosak, Vladimir  S. Matveev\thanks{ Institute of Mathematics, FSU Jena, 07737 Jena Germany,  vladimir.matveev@uni-jena.de}}
\begin{document}
\maketitle

\begin{abstract} 
We show that if two 4-dimensional  metrics  of arbitrary signature on one manifold  are geodesically equivalent (i.e., have the same geodesics considered as unparameterized curves) and are solutions of the Einstein field equation with the same stress-energy tensor, then they are affinely equivalent or flat.     Under the additional assumption  that 
 the metrics are complete or the manifold  is closed, the result survives in all dimensions $\ge 3$. 
\end{abstract} 

\section{Definitions and results}
\label{results} 

Let $(M^n,g)$ be a connected  pseudo-Riemannian  manifold of arbitrary signature  of dimension $n\ge 3$.  

 We say that a metric $\bar g$ on $M^n$ is \emph{geodesically equivalent} to $g$, if every geodesic of $g$ is a (possibly, reparametrized) geodesic of $\bar g$.  We say that they are \emph{affinely  equivalent},  if the Levi-Civita connections of $g$ and $\bar g$  coincide.

 In this paper we study the question whether two geodesically equivalent metrics 
  $g$ and $\bar g$  can satisfy  
the  Einstein field equation with the same  stress-energy tensor:
\begin{equation} \label{stress} 
R_{ij}-  \tfrac{R}{2} \cdot g_{ij}=\bar R_{ij}-  \tfrac{\bar R}{2} \cdot \bar g_{ij},
\end{equation} 
 where $R_{ij}$ ($\bar  R_{ij}$, respectively) 
  is the Ricci tensor of the metric $g$ ($\bar g$, respectively), 
  and $R:= R_{ij}g^{ij}$ ($\bar R:= \bar R_{ij}\bar g^{ij}$, respectively, $\bar {g}^{k \ell}$ is  the tensor dual to $\bar g_{i j}$:  \ $\bar {g}^{s i}\bar g_{s j}= \delta_{\ j}^i$)  is the scalar curvature. 

There exist the following trivial examples of such a situation: 

\begin{enumerate}\item If geodesically equivalent metrics 
$g$ and $\bar g$ are flat, then their stress-energy tensors vanish identically and therefore coincide. Examples of geodesically equivalent flat metrics are classically known and  can be constructed as follows: take the classical projective transformation $p$  of 
$(U\subseteq  \R^n, g_{\textrm{\tiny standard}})$ (i.e., a local diffeomorphism that takes straight lines to straight lines,  there is a $(n^2+2n)$-dimensional group of it) and   consider 
the pullback of the standard euclidean metric  $g_{\textrm{\tiny standard}} $; $\bar g = p^*g_{\textrm{\tiny standard}}$.  
  It is clearly flat and      geodesically equivalent to the initial metric $g_{\textrm{\tiny standard}}$. If $p$ is not a classical affine transformation (the subgroup of affine transformations is $n^2 + n$-dimensional), $\bar g$ is  not affinely equivalent to $g_{\textrm{\tiny standard}}.$

\item If $g$ and $\bar g$ are affinely equivalent metrics with vanishing scalar curvature, then their stress-energy  tensors coincide with the Ricci tensors and therefore coincide 
(since even Riemannian curvature tensors coincide). There are many examples of such a situation, a possibly simplest one is as follows:  
Take an arbitrary metric $h= h_{ij}, \ i,j= 2,...,n$ of zero scalar curvature on $\R^{n-1}(x_2, ..., x_{n-1})$ and consider the direct product metric  $g= dx_1^2 + \sum_{i,j=2}^n h_{ij}dx_idx_j$ on  $\R^n = \R(x_1)\times \R^{n-1} (x_2,...,x_n)$.  Then, for this 4-dimensional metric, and also for the (affinely equivelent) metric $g= dx_1^2 + 2 \sum_{i,j=2}^n h_{ij}dx_idx_j$,  the scalar curvature  is zero.

\item The metric  $\bar g:= \const\  \cdot g$ has the same  stress-energy tensor as $g$. Indeed, $R_{ij}= \bar R_{ij}$, and $\bar R:= \bar g^{ij} R_{ij}= \tfrac{1}{\const} R$ so $R g_{ij}= \tfrac{1}{\const} R \cdot \const\ g_{ij}= \bar R \bar g_{ij}$.    
\end{enumerate}

In the present paper we show that in dimensions  3 and 4  this list of trivial examples contains all possibilities: 
\begin{Th} \label{main}
 If two geodesically equivalent metrics  $g$ and  $\bar g$ on a connected $M$  of dimension 3 or 4    satisfy  \eqref{stress}, then one of the following possibilities takes place:
  \begin{enumerate} 
  \item $g$ and $\bar g$ are affinely  equivalent metrics with zero scalar curvature, or    
  \item $g$ and $\bar g$ are flat,  or   
   \item $\bar g= \const\ g$ for a certain $ \const \in \mathbb{R}$   
  \end{enumerate}

\end{Th}

By this theorem, 
 unparameterized geodesics determine the Levi-Civita  connection of a  3 or 4-dimensional   metric uniquely within the solutions of the Einstein field equation with the same stress-energy tensor provided the metric  is  not flat. 

The motivation to study this question came from physics. It is known that 
 geodesics of a space-time metric   correspond to the trajectories of the free falling uncharged particles, and that certain astronomical observations give  the trajectories  of  free falling uncharged particles as   unparameterized curves; moreover, unparameterized geodesics   and how and whether they determine  the metric   were actively studied by theoretical physicists   (cf \cite{EPS72,Petrov1,veblen23,Weyl1})
   in the context of general relativity.  The space-time metric is a solution  of the   Einstein equation (there of course could be many solutions of the Einstein equation with the same stress-energy tensor)
  and  our theorem implies that if we know the (unparameterized) trajectories   of  free falling uncharged particles and the stress-energy tensor, then we know (i.e., can in theory reconstruct) the metric or at least the Levi-Civita connection of the metric. 
 
 The dimension 4 is probably the dimension that could be interesting for  physics, since  space-time metrics are naturally   4-dimensional. The result for dimension 3 is essentially easier; that's why we put it here.  
In dimension two, the stress-energy  tensor of every metric  is identically zero and (the analog of) Theorem \ref{main}  is evidently wrong.  
It is also  wrong in higher dimensions, we show a  example in dimensions  $\ge 5$. 
The metrics $g$ and $\bar g$ in this example both have zero scalar curvature and their Riemannian curvature tensors coincide. We do not know whether all geodesically equivalent not affinely equivalent  metrics with the same stress-energy tensors have zero scalar curvature, but can show that the scalar curvature must be constant.

\begin{Th} \label{main2}
Suppose two nonproportional   geodesically equivalent metrics  $g$ and  $\bar g$ on a connected $M^n$, $n\ge 5$
   satisfy  \eqref{stress}. Then, the scalar curvatures of the metrics are  constant. 
 \end{Th} 

Combining this theorem  with \cite{KM,Mounoud2010} we obtain that  
in the global setting, when the manifold is closed (= compact without boundary), or  when both metrics are complete, the analog of Theorem \ref{main} is still true in all dimensions.

We say that a (complete in both directions)  $g$-geodesic   $\gamma:\R\to M$ is $\bar g$-complete, of there exists a diffeomorphism $\tau:\R\to \R$ such that the curve $\bar \gamma:= \gamma\circ \tau$ is a $\bar g$-geodesic.

\begin{Cor} \label{cor1}  Suppose  geodesically equivalent  metrics  $g$ and  $\bar g$ on a  connected $M^n$, $n\ge 5$, such that $g$ has indefinite signature 
  satisfy  \eqref{stress}. 
 Assume in addition that   every light-like     $g$-geodesic  $\gamma$   is complete in both direction   and is  $\bar g$-complete.  
  Then, the metrics are affinely equivalent.
\end{Cor} 

 \begin{Cor} \label{cor2}  Suppose  two  geodesically complete  geodesically equivalent  metrics  $g$ and  $\bar g$ on a  connected $M^n$, $n\ge 5$, such that $g$ is positively definite or negatively definite, 
  satisfy  \eqref{stress}. 
  Then, the metrics are affinely equivalent.
\end{Cor} 

\begin{Cor} \label{cor3} Suppose  two geodesically equivalent    metrics  $g$ and  $\bar g$ on a  closed   connected $M^n$, $n\ge 5$, 
  satisfy  \eqref{stress}. 
  Then, the metrics are affinely equivalent.\end{Cor}

  Probably the most famous special case of Theorem \ref{main} that was known before is due to A. Z. Petrov  \cite{Petrov1} (see also \cite{Hall2007} and \cite{einstein}): he has shown that    
      {\it   4-dimensional Ricci-flat nonflat  metrics  of Lorentz signature can not be geodesically  equivalent, unless they are affinely  equivalent. }
       It is one of the results  Petrov  obtained in 1972 the Lenin   prize, 
       the most important  scientific award of the  Soviet Union,   for. 

\section{Proof of Theorems 1 and 2.}
\subsection{Plan of the proof.}

We  start with  recalling in  \S \ref{standard}  certain  known facts from the theory of geodesically equivalent metrics that will be used in the proof. In \S \ref{sec2}  we prove an important technical statement: we show that if the minimal polynomial of the  tensor $a^i_{\ j}$ defined by \eqref{a} has degree 2, then   geodesically equivalent metrics that were used to construct $a^i_{\ j}$ are  warped product  metrics provided they are not affinely equivalent. In \S  \ref{sec3} we prove Theorem \ref{main} for geodesically equivalent warped product metrics.  

In \S \ref{vnb} we use the connections between the Ricci tensors of geodesically equivalent metrics 
 to derive the formula \eqref{vb} which will play an important role  in the proof.

\weg{By \S \ref{vnb}, $\lambda_{i,j}= \mu g_{ij}+ \frac{R}{2(n-1)} a_{ij}$.}

The proof  depends on the behavior of the scalar curvature of a metric:  the following three cases use different ideas:   

\begin{itemize}\item Case 1: $R=\const \ne 0$. 
\item Case 2: $dR\ne 0$.
\item Case 3: $R=0$.

\end{itemize}

Clearly,  almost every point of $M^n$ belongs to one of the cases 1,2,3; so it is sufficient to prove Theorem \ref{main} 
under the assumption of these cases. We will do it in \S\S \ref{case1}, \ref{case2}, \ref{case3} respectively.   
The first and the second    cases will be reduced to the warped product case solved in \S \ref{sec3}, but in each case the reduction will be different. In the second case, and also in the ``warped product part" (i.e., in \S    \ref{sec3}) we will  work in arbitrary  dimensions $n \ge 3$, so we simultaneously prove  Theorem \ref{main2}.

 \subsection{ Standard formulas we will use} \label{standard} 
 We work in tensor notations with the background metric $g$. That means, we sum with respect to repeating indexes, use $g$ for  raising and lowing  indexes (unless we explicitly mention), and use the Levi-Civita connection of $g$ for   covariant differentiation which we denote by comma.

 As it was known already to Levi-Civita \cite{Levi-Civita},  two connections $\Gamma= \Gamma_{jk}^i $ and  $\bar \Gamma= \bar \Gamma_{jk}^i $  have the same unparameterized geodesics, if and only if their difference is a pure trace: there exists a $(0,1)$-tensor $\phi_i  $ 
such    that \begin{equation} \label{c1} 
 \bar \Gamma_{jk}^i  = \Gamma_{jk}^i + \delta_{\ k}^i\phi_{j} + \delta_{\ j}^i\phi_{k}.    
   \end{equation}

If $\Gamma$  and   $\bar \Gamma$  related by  \eqref{c1}  are Levi-Cevita connections of  metrics $g$ and $\bar g$, then one can find explicitly (following Levi-Civita \cite{Levi-Civita}) a function $\phi$ on the manifold such that its differential $\phi_{,i}$  coincides with the $(0,1)$-tensor  $\phi_i$: indeed, contracting \eqref{c1}  with respect to  $i$ and $j$, we obtain 
 $\bar \Gamma_{s i}^s   = \Gamma_{s i}^s  + (n+1) \phi_{i}$. 
  From the  other side, for the Levi-Civita 
 connection  $\Gamma$ of a metric $g$  we have  $  \Gamma_{s k}^s  = \tfrac{1}{2} \frac{\partial \log(|det(g)|)}{\partial x_k} $.  Thus, 
 \begin{equation} \label{c1,5}  \phi_{i}= \frac{1}{2(n+1)}  \frac{\partial }{\partial x_i}  \log\left(\left|\frac{\det(\bar g)}{\det( g)}\right|\right)= \phi_{,i} \end{equation} 
  for the function $\phi:M\to \mathbb{R}$ given by 
  \begin{equation} \label{phi}  \phi:= \frac{1}{2(n+1)} \log\left(\left|\frac{\det(\bar g)}{\det( g)}\right|\right). \end{equation}  In particular, the derivative of $\phi_i$ is  symmetric, i.e., $\phi_{i,j}= \phi_{j,i}$.

The formula  \eqref{c1}   implies  that two metrics $g$ and $\bar g$ are geodesically equivalent if and only if   for a certain $\phi_{i}$ (which is, as we explained above, the differential of $\phi$ given by  \eqref{phi}) we have 
\begin{equation}\label{LC}
    \bar g_{ij, k} -  2 \bar g_{ij} \phi_{k}-  \bar g_{ik}\phi_{j} -   \bar g_{jk}\phi_{i}= 0, 
\end{equation} 
where ``comma" denotes the covariant derivative with respect to the connection $\Gamma$. 
Indeed, the left-hand side of this equation is the covariant derivative with respect to  $\bar \Gamma$, and vanishes if and only if  $\bar \Gamma$ is the Levi-Civita connection for $\bar g$.   
Clearly, the metrics $g$ and $\bar g$ are affinely equivalent, if $\phi_i\equiv 0$, or, which is the same, if $\phi= \const$.

The equations \eqref{LC} should be viewed a system of  PDE on the unknowns $\bar g_{ij}$ and $\phi_i$. It  can be linearized by a clever  substitution (which was already known to R.  Liouville \cite{liouville} and Dini \cite{Dini}  in dimension 2 and is due to Sinjukov \cite{sinjukov} is other dimensions, see also \cite{benenti,eastwood}):  consider     $a_{ij}$ and $\lambda_i$ given by 
\begin{eqnarray} \label{a}
a_{ij} &=   &e^{2\phi} \bar g^{s q} g_{s i} g_{q j},\\  \label{lambda} 
\lambda_{i} & = &  -e^{2\phi}\phi_s \bar g^{s p} g_{p i}, \end{eqnarray}
where $\bar {g}^{s p}$ is  the tensor dual to $\bar g_{i j}$:  \ $\bar {g}^{s i}\bar g_{s j}= \delta_{\ j}^i$. 
It is an easy exercise to show that the following linear  equation on  the symmetric $(0,2)-$tensor $a_{ij}$ and $(0,1)-$tensor $\lambda_i$ is     equivalent to \eqref{LC} 
 \begin{equation} \label{basic} 
 a_{ij,k}= \lambda_i g_{jk} + \lambda_j  g_{ik}. 
 \end{equation}

Note that there exists  a  function $\lambda$ such that its differential is precisely 
 the $(0,1)-$tensor $\lambda_i$: indeed, multiplying \eqref{basic}  by $g^{ij}$ and summing with respect to repeating indexes $i,j$ we obtain $(g^{ij}a_{ij})_{,k} = 2  \lambda_k$. Thus,
$\lambda_i$ is the differential of the function 
\begin{equation}\label{lam} 
\lambda:= \tfrac{1}{2} g^{qp }a_{qp}.   
\end{equation}  
In particular, the covariant derivative of $\lambda_i$ is symmetric:  
 $\lambda_{i,j} = \lambda_{j,i}$. Clearly, the metrics $g$ and $\bar g$ are affinely equivalent, if $\lambda_i\equiv 0$, or, which is the same, if $\lambda= \const$.

\begin{Rem}  In this paper an important role plays the tensor $A:= a^{i}_{ \ j}$ which we will view as a field of endomorphims of $TM$; combining the formulas \eqref{a} and \eqref{phi} we see that it is given by the formula 
\begin{equation}
\label{L}
A= a_{\ j}^i := \left(\left|\frac{\det(\bar g)}{\det(g)}\right|\right)^{\frac{1}{n+1}} \bar g^{is}
 g_{sj}.
\end{equation}
 One can reconstruct (up to the sign but since the equation \eqref{stress} survives if we replace $\bar g$ by $-\bar g$, the sign in not essential)
 the metric $\bar g$ (considered as a bilinear form) by $A$ and $g$ by the formula 
\begin{equation}
\label{LG}
\bar g(  \ , \ ) = \frac{1}{|\det(A)|}g( A^{-1} \ , \ ).
\end{equation}
\end{Rem}

 Integrability conditions for the equation \eqref{basic}  (we substitute the derivatives of $a_{ij}$ given by \eqref{basic} in the formula  $a_{ij,\ell k}- a_{ij,k\ell}= a_{i s }R^{s}_{\ jk\ell} +  a_{s j}R^{s}_{\ ik\ell}$, which is true for every $(0,2)-$tensor  $a_{ij}$)   were  first obtained by Solodovnikov \cite{s1} and are   
 
\begin{equation}  a_{i s }R^{s}_{\ jk\ell} +  a_{s j}R^{s}_{\ ik\ell} =\lambda_{ \ell,i} g_{jk}+\lambda_{ \ell,j} g_{ik}-\lambda_{ k,i} g_{j\ell}-\lambda_{ k,j} g_{i\ell}. \label{int1} \end{equation}

For further use let us recall the  following well-known fact which can also be obtained  by simple calculations (the straight-forward way is to  replace $\Gamma$ by $\bar \Gamma$ given by  \eqref{c1} in the formula for the Riemannian  curvature 
and then for the  Ricci tensor):  the 
Ricci-tensors of connections related by \eqref{c1} are  connected by the formula 
\begin{equation} \label{ric} \bar R_{ij} = R_{ij}- (n-1)(\phi_{i,j}- \phi_{i}\phi_{j}),  
\end{equation} 
 where $R_{ij}$ is the Ricci-tensor  of $\Gamma  $ and $\bar R_{ij}$ is the Ricci-tensor  of $\bar \Gamma $.

 Important special case of the metrics we will consider in our proof  will be the metrics such that they admit a solution $(a_{ij}, \lambda_i)$ of \eqref{basic} with $a_{ij}\ne \const \cdot g_{ij}$, such that the derivative  of  $\lambda_i$ satisfies, for a certain constant $B$ and for a certain function $\mu$, the equation 
\begin{equation} \label{vn}
 \lambda_{i,j}= \mu g_{ij} + B a_{ij}.
 \end{equation}
 This condition may  look artificial from the first glance, but it is  not, since it naturally appears in many situations in the theory of geodesically equivalent metrics. For example, if $g$ is Einstein, then every solution $(a_{ij}, \lambda_i)$ 
 satisfies this condition (with $B= - \frac{R}{n(n-1)})$, see \cite[Eq. (24)]{einstein}. Moreover, if the dimension  of the space of solutions of \eqref{basic} is at least three, then there exists a constant $B$ such that every solution of \eqref{basic} satisfies \eqref{vn} (the constant $B$ is the same for all solutions but the function $\mu$ depends on the solution), see  \cite[Lemma 3]{KM}. Moreover, the constant $B$ is unique for all solutions and is the same on the whole (connected) manifold \cite[\S\S 2.3.4, 2.3.5]{KM}.    In our setting, under the assumption that the scalar curvature $R$ is a constant, the equation \eqref{stress} implies the equation \eqref{vn} for the constant $B= - \tfrac{R}{2(n-1)}$, see \S \ref{vnb}.

Moreover, if \eqref{vn} is satisfies, then 
the function $\mu$ necessary satisfies  the equation $\mu_{,i}=2B\lambda_i$ (see \cite[Rem. 10]{KM}), so the triple $(a, \lambda, \mu)$ satisfies the following Frobenius-type system:   
  \begin{gather}
    \label{eq:mg+Ba}
    \left\{
      \begin{array}{ll}
        a_{i j,k}&=\lambda_i g_{j k} + \lambda_j g_{i k}\\
        \lambda_{i,j}&=\mu g_{i j}+B a_{i j}\\
        \mu_{,i}&=2B\lambda_i
      \end{array}
    \right. 
  \end{gather}

For further use we need the following 

 \begin{Lemma}[cf Lemma 9 of \cite{KM}] \label{const}
 Let $g, \bar g$ be geodesically equivalent  metrics on a connected $M^{n\ge 3}$. Assume that  the metric $g$   admits a solution  $(a_{ij}, \lambda_i)$  with $\lambda_i\ne 0$ of \eqref{basic} 
 such that \eqref{vn} holds. Assume also that  the metric $\bar g$ admits a solution  $(\bar a_{ij}, \bar \lambda_i)$ of the natural analog of \eqref{basic} with 
 $\bar \lambda_i\ne 0$   such that the natural analog of \eqref{vn} holds; we denote the natural analog of $B$ by $\bar B$.

   Then, 
   the following formula holds:
 \begin{equation}\label{f1} 
 \phi_{i,j} - \phi_{i} \phi_{j} = -B g_{ij}+\bar B \bar g_{ij}  .
 \end{equation} 
 \end{Lemma}
 
 {\bf Proof. }   
 We covariantly  differentiate \eqref{lambda} (the index of differentiation is ``j"); then we substitute the expression \eqref{LC} for $\bar g_{ij,k}$   to obtain 
 \begin{equation} \label{f21} \begin{array}{ccl}
 \lambda_{i,j} &=& -2 e^{2\phi}\phi_{j} \phi_p \bar g^{p q} g_{q i}-e^{2\phi}\phi_{p,j} \bar g^{p q} g_{q i}+e^{2\phi}\phi_p  \bar g^{p s} \bar g_{s l,j} \bar g^{l q} g_{q i} \\ &=&   -e^{2\phi}\phi_{p,j} \bar g^{p q} g_{q i}+e^{2\phi}\phi_p \phi_s \bar g^{p s}   g_{ i j }+  e^{2\phi}\phi_{j} \phi_l \bar g^{lq}g_{q i}      \end{array}    ,  
 \end{equation} 
 where $\bar g ^{p q}$ is the tensor dual to $\bar g_{p q}$, i.e.,  $\bar g ^{p i}\bar g_{p j}= \delta_{\ j}^i$ . 
 We now  substitute  $\lambda_{i,j}$ from \eqref{vn},  use that $a_{ij}$ is given by \eqref{a}, and divide by $e^{2\phi}$ for cosmetic reasons   to  obtain 
 \begin{equation} \label{f31} 
 e^{-2\phi} \mu g_{ij} + B \bar g^{p q} g_{p j}g_{q i} = -\phi_{p, j} \bar g^{p q} g_{q i}+\phi_p \phi_{s} \bar g^{p s} \bar  g_{ i j }+ \phi_{j} \phi_{l} \bar g^{lq}g_{q i}.  
 \end{equation}  
 Multiplying with $g^{i m } \bar g_{m k}$,  we obtain 
 \begin{equation} \label{f4} 
 \phi_{k,j}-\phi_{k}\phi_{j} =\underbrace{(\phi_p \phi_q \bar g^{p q} - e^{-2 \phi } \mu )}_{\bar b} \bar g_{kj} - B g_{kj}.  
  \end{equation} 
   The same holds with the roles of $g$ and $\bar g$ exchanged  (the function \eqref{phi} constructed by  the interchanged pair $\bar g, g$ is evidently equal to $-\phi$).   We  obtain 
  \begin{equation} \label{f4b}
   -\phi_{k;j}-\phi_{k}\phi_{j} =\underbrace{(\phi_p \phi_q g^{p q} - e^{2 \phi } \bar \mu )}_{b}  g_{kj} - \bar B \bar g_{kj},
  \end{equation}  
  where $\phi_{i;j}$ denotes the covariant derivative  of $\phi_i$ with respect to the Levi-Civita connection of the metric $\bar g$. 
Since the Levi-Civita connections of $g$ and of $\bar g$ are related by the formula  
(\ref{c1}), we have 
$$
  -\phi_{k;j}-\phi_{k}\phi_{j} = \underbrace{-\phi_{k,j} + 2 \phi_k\phi_j}_{ -\phi_{k;j}} - \phi_{k}\phi_{j}= -(\phi_{k,j}-\phi_{k}\phi_{j}).$$
  We see that the left hand side of \eqref{f4} is equal to minus the left hand side of \eqref{f4b}.
  Thus, $b\cdot g_{ij} - \bar B \cdot \bar g_{ij} = B\cdot g_{ij} - \bar b \cdot \bar g_{ij}$ holds. Since the metrics $g$ and $\bar g$ are not proportional  by assumption, $\bar b= \bar B$ as we explained above, and the formula \eqref{f4} coincides with \eqref{f1}.  Lemma is proved.

\begin{Rem} \label{remnew}
We see that under the assumptions of Lemma \ref{const} the constant $B$ is given in view of \eqref{f4b} by  $$B=\phi_p \phi_q g^{p q} - e^{2 \phi } \bar \mu .$$
\end{Rem}

   \subsection{Geodesically equivalent metrics such that the minimal polynomial of  $A= a^i_{\ j}$ has degree 2.} \label{sec2}

   Assume that $(a_{ij},\lambda_i)$ is a nontrivial (i.e., $\lambda_i\ne 0$) solution of \eqref{basic}. We assume $n=\dim(M)\ge 3$. We will discuss the situation when 
    the minimal polynomial of the $(1,1)$-tensor $A= a^i_{\ j}$ (viewed as an endomorphism of $TM$)  has degree  at most 2 (in every  point of some neighborhood), i.e., when there exist functions $c_1$ and $c_2$ such that 
    \begin{equation}\label{minimal}  A^2 + c_1\ A + c_2\ \Id = a^i_{\ k} a^k_{ \ j} + c_1 a^i_{\ j} + c_2 \delta^i_{\ j} =0.\end{equation}  In other words, we assume that $A$ has the following real Jordan normal form (at every point of the neighborhood we are working in); in all matrices below we assume that zeros stay on the empty spaces and all diagonal blocks are square matrices

    \begin{equation} \label{matrices1}
      \left(\begin{array}{ccc|ccc}
         &&&&&\\
         & \rho_1\Id_{k\times k} &  &  & &\\
         &&&&&\\
        \hline
        &&&&& \\
        &&&&\rho_2\Id_{(n-k)\times (n-k)}&  \\
        &&&&&
      \end{array}\right), \ \end{equation}  \begin{equation} \label{matrices2}  \left(\begin{array}{cc|ccc}
         \rho & 1&&&\\
         & \rho &    & &\\
        \hline
        &&&& \\
        &&&\rho\Id_{(n-2)\times (n-2)}&  \\
        &&&&
      \end{array}\right),\   \left(\begin{array}{cc|c|cc}
         \alpha& \beta &&&\\
         -\beta & \alpha   & & &\\ \hline
         &&\ddots &&\\
        \hline
       &&& \alpha& \beta \\
     & & &   -\beta & \alpha   
      \end{array}\right),
     \end{equation}
       where $\Id_{k\times k}= \textrm{diag(1,...,1)}$ denotes the matrix of the identity endomorphism of $\mathbb{R}^k; $
 we assume $0<k<n$.       
    
   Our goal  is to prove the following 
    \begin{Lemma}  \label{lem1}  Let $(M^n, g)$ be a pseudo-Riemannian connected manifold of dimension $n\ge 3$ and $(a_{ij}, \lambda_i)$ be a solution of \eqref{basic}  such that  $\lambda_i\ne 0$ 
      such that   there exist functions  $c_1, c_2$ on $M$ such that \eqref{minimal} holds. Then, at the generic  point of $M$  the Jordan normal form of $A= a^i_{\ j}$ is as in \eqref{matrices1}, moreover $k= 1$ or $k =n-1$; in other words, the Jordan normal form of $A= a^i_{\ j}$ has no Jordan blocks of dimension $\ge 2$, all eigenvalues of A= $a^i_{\ j}$ are real, and at the points where there are two eigenvalues,  one of them has algebraic multiplicity one.  Moreover, the other eigenvalue of  $A$  (considered as a function on $M$) is constant. 

Moreover,          in a neighborhood of 
       every  point  such that $a_{ij}$ is not proportional to $   g_{ij}$,    there exists a coordinate system $(x_1,...,x_n)$  where  the matrices  of $g$  and of $a^i_{\ j}$ are given by 
       
       \begin{equation} \label{wp} 
      g_{ij}= \left(\begin{array}{c|ccc}
          \pm 1  &&&\\
        \hline
        &&& \\
        &&\sigma(x_1)h &  \\
        &&&
      \end{array}\right), \ 
      a_{\ j}^i=  \left(\begin{array}{c|ccc}
          \sigma(x_1) +C  &&&\\
        \hline
        &&& \\
        &&C\ \Id &  \\
        &&&
      \end{array}\right),
       \end{equation}
where $\sigma$ is a function of $x_1$, $C$ is a constant, and $h$ is a symmetric nondegenerate $(n-1)\times (n-1)$-matrix whose entries depend on $x_2,...,x_n$.       \end{Lemma} 
    
    {\bf Proof.}  The proof is based on the Splitting and Gluing Lemmas from \cite{splitting} and on \cite[Proposition 1]{local}.  
    If the minimal polynomial of $A$ has degree $1$ in a neighborhood of a point, 
    the metric $g$ is conformally equivalent  to $\bar g$;  by the classical result  of Weyl \cite{Weyl2} the conformal coefficient is a constant. Then, by   \cite[Proposition 1]{local}, the metrics are proportional with constant coefficient on the whole manifold (assumed connected)  which contradicts the assumptions. 
      
    Assume the Jordan form of $A=a^i_{\ j}$ is as the first  one in \eqref{matrices2} in a small neighborhood. Then, the geometric multiplicity  (i.e., the dimension of the eigenspace) of the eigenvalue $\rho$ is $\ge 2$  implying by \cite[Proposition 1]{local} that (the function)  $\rho$ is constant  on the whole manifold  so $\lambda= \tfrac{1}{2}a^s_{ \ s}=  \tfrac{n}{2} \rho$ is also constant so $\lambda_i=0$. 
    Similarly, if the Jordan form of $a^i_{\ j}$ is as the second  one in \eqref{matrices2} in a small neighborhood, the geometric multiplicity of the eigenvalues $\alpha \pm i \beta$ is $\ge 2$ implying by 
    \cite[Proposition 1]{local} that (both functions)  $\alpha$ and $\beta$   are  constant   so $\lambda= \tfrac{1}{2} a^s_{\ s }=  \tfrac{n}{2} \alpha$ is also constant so $\lambda_i=0$. 
Assume now the Jordan normal form of $A= a^i_{\ j}$ is as in \eqref{matrices1}. If $k\ne 1$ and  $k\ne n-1$, then the geometric 
multiplicities of both eigenvalues  are $\ge 2$ so we again obtain that $\lambda= \tfrac{1}{2} a^i_{i}=  \tfrac{k}{2} \rho_1 +  \tfrac{n-k}{2} \rho_2$ is constant so $\lambda_i=0$. Thus, $k= 1$ or $n-k= 1$; without loss of generality we assume $k=1$; then $\rho_2=  \const$.  

The characteristic polynomial of $A= a^i_{\ j}$ is $\chi= (t- \rho_1) (t- \rho_2)^{n-1}$; we denote $(t-\rho_1) $  by $\chi_1$ and $(t- \rho_2)^{n-1}$ by $\chi_2$. The factorization $\chi= \chi_1\cdot \chi_2$ is  admissible in the terminology of \cite[\S 1.1]{splitting}. Then, by the Splitting  \cite[Theorem  3]{splitting} and  Gluing Lemmas \cite[Theorem  4]{splitting},  there exists  
a  coordinate system $(x_1,...,x_n)$  such that in this coordinates the eigenvalue $\rho_1$ is a function of $x_1$  and   the matrices of $a^i_{\ j}$ and of $g$ are as follows:
\begin{equation} \label{a1} 
      a^{i}_ {\ j} =  \left(\begin{array}{c|ccc}
          \rho_1  &&&\\
        \hline
        &&& \\
        &&\rho_2 \Id &  \\
        &&&
      \end{array}\right), \  \ \end{equation}\begin{equation} \label{a1bis}   
  g_{ij}=      \left(\begin{array}{c|ccc}
          f(x_1) \chi_2(\rho_1)  &&&\\
        \hline
        &&& \\
        &&\chi_1(\rho_2) h &  \\
        &&&
      \end{array}\right)=  \left(\begin{array}{c|ccc}
          f(x_1) (\rho_1- \rho_2)^{n-1}  &&&\\
        \hline
        &&& \\
        &&(\rho_2- \rho_1) h &  \\
        &&&
      \end{array}\right) , \ 
\end{equation}
where $h$ is a symmetric nondegenerate $(n-1)\times (n-1)$ matrix whose components depend on the variables $x_2,...,x_n$ only and $f$ is a function of $x_1$. Replacing the first coordinate by an appropriate function $X_1$ of it  (such that 
$dX_1=  \sqrt{|f(x_1) (\rho_1- \rho_2)^{n-1} |} dx_1$)    we can make the (1,1)-component of $g$ to be $\pm 1$. Finally   we see 
 that $g$ and $a^i_{\ j}$ are given by the  formulas \eqref{wp} 
  with $ \sigma=  \rho_1$ and  $C= \rho_2$.  Lemma is proved. 

\begin{Rem} Assume in addition that  the solution $a_{ij}$ came from a metric $\bar g_{ij}$ by  \eqref{a}, i.e., assume that  $a_{ij}$ in nondegenerate, i.e., assume that $\rho_2\ne 0$ and $\rho_1(x_1)\ne 0$ at every point of the neighborhood  we are working in. Then,    by \eqref{LG}, the 
matrix of $\bar g$ is given by 

\begin{equation} \label{barwp} \bar g_{ij}= \frac{1}{C^{n-1}}   \left(\begin{array}{c|ccc}
          \pm \tfrac{ 1 }{ (\sigma +C)^2}  &&&\\
        \hline
        &&& \\
        &&\frac{\sigma}{(\sigma +C)C} h &  \\
        &&&
      \end{array}\right) . \ 
\end{equation}

Without loss of generality we can assume later then the sign $\pm$ of the $(1,1)$-entry  of $g$ and of $\bar g$ is ``$+$'', since multiplication of  $g$ or  of $\bar g$  by $-1$ does not affect the equation \eqref{stress}.

\end{Rem}

    \subsection{ Proof for geodesically equivalent warped product metrics. } \label{sec3}
    
    We will now prove Theorem \ref{main} under the additional  assumption that the geodesically equivalent metrics $g$ and $\bar g$ satisfying \eqref{stress} are given by the formulas 
    (\ref{wp}, \ref{barwp}). We prove 
    
    \begin{Lemma} \label{lem2}
    Assume the metrics  $g$ and $\bar g$ given by (\ref{wp}, \ref{barwp}) satisfy \eqref{stress}. Then, the function $\sigma$ is a constant, so the metrics are affinely equivalent.  
    \end{Lemma} 
    
    {\bf Proof. } We prove the Lemma by direct calculations: a straightforward way to do it (at least in the 3- and 4-dimensional case which will be used  in the proof of Theorem \ref{main})   is to use any computer algebra program, for example Maple, to calculate the difference between the left- and the right hand sides of \eqref{stress}. One immediately sees that the $i,j$-component of the difference with $i\ge 2,  j \ge 2$ is proportional to  the corresponding entry of $h_{ij}$ with the same coefficient of  the proportionality which is proportional 
    to  $(\sigma')^2$.  Since it is zero by assumptions, $\sigma$ is a constant and the metrics are affinely equivalent. 
  
As a part of the proof of Theorem \ref{main2} we need this calculation in arbitrary dimension;    
let us explain  a small trick that  helps to  calculate the difference between the left- and the right hand sides of \eqref{stress} `by hands' and in any dimension. 
   
   We will use that the conformally equivalent metric $\frac{1}{\sigma} g $ is the direct product metric so  its Ricci tensor has the form 
 \begin{equation}  \left(\begin{array}{c|ccc}
          0 &&&\\
        \hline
        &&& \\
        &&H_{ij} &  \\
        &&&
      \end{array}\right)\label{curvH} \end{equation}
   where  $H$ is the Ricci-tensor of the $(n-1)$-dimensional metric $h_{ij}$ (viewed as a metric on $U\subseteq \R^{n-1}(x_2,...,x_n)$, and its scalar curvature is simply the scalar curvature of $h_{ij}$.  
   Now, it is well known that the Ricci-tensors and the scalar curvatures  of any the conformally equivalent metrics  $g $ and 
   $ \hat g:= e^{-2\psi} g$ 
   are related by 
   \begin{equation}\label{conf} 
   \begin{array}{ccl}  \hat R_{ij} &=&  R_{ij} - (n-2)(\psi_{i,j}- \psi_i \psi_j) - (\Delta_2  + (n-2)\Delta_1)g_{ij} , \\ \hat R&=& -e ^{-2\psi} ( R+ 2(n-1)\Delta_2 + (n-1)(n-2)\Delta_1 ), 
  \end{array} \end{equation}
   where $\Delta_2$ is the Laplacian of $\psi$, $\Delta_2= \psi_{i,j} g^{ij}$, and $\Delta_1$ is the square of the  length of $\psi_i$ in $g$, $\Delta_1:= g^{ij}\psi_{,i} \psi_{,j}$. 
   In our case the  role of the metric $g$ in \eqref{conf}  plays the direct product metric $\frac{1}{\sigma} g $ and $\phi= -\tfrac{1}{2} \log |\sigma|.$ After some relatively simple calculations   we obtain  $R_{ij}$ as  an algebraic  expression in  $H_{ij}$, $h_{ij}$, $\sigma$, $\sigma'$ and $\sigma''$,    and also  $R $ as an algebraic  expression in  $H= H_{ij} h^{ij}$,  $\sigma$, $\sigma'$ and $\sigma''$.

   Similarly, the metric  $\frac{C(\sigma+ C)} {\sigma}\bar  g $  which is conformally equivalent to the metric $\bar g$ 
   is also the direct product metric so its Ricci curvature also is  as in \eqref{curvH}. We again combine it with  \eqref{conf} and calculate the scalar and the Ricci curvatures of $\bar g$.   Substituting 
    the result of the calculation in
 the left hand side of    \eqref{stress} minus the right-hand side of \eqref{stress}, and considering the components of the result for $i,j\ge 2$, we see that $H$ and $H_{ij}$   disappear and we obtain the following  condition on the function $\sigma$ only
$$
\frac{(n-2)(n-1) (\sigma')^2}{6\sigma(\sigma + C)}=0.
$$  
   Then,  $\sigma'= 0$, which implies that $
  \sigma$ is a constant and the metrics are affinely equivalent.

\subsection{ $\lambda_{i,j}$ is a linear combination of $g_{ij}$ and $a_{ij}$ (with functional coefficients)} \label{vnb}
Assume geodesically equivalent $g$ and $\bar g$ on $M^n$ satisfy \eqref{stress}. 
Rearranging the terms in  \eqref{stress}, we obtain 
  $R_{ij} - \bar R_{ij} = \tfrac{R}{2}g_{ij} - \tfrac{\bar R}{2}g_{ij}$. 
Substituting \eqref{ric} inside, we obtain
\begin{equation}\label{f11} 
 \phi_{i,j} - \phi_{i} \phi_{j} = \tfrac{R}{2(n-1)} g_{ij}-\tfrac{\bar R}{2(n-1)} \bar g_{ij}  ,
 \end{equation} 

Now we  covariantly  differentiate \eqref{lambda} (the index of differentiation is ``j"); then we substitute the expression \eqref{LC} for $\bar g_{ij,k}$,  and finally we substitute \eqref{f11}   to obtain 
 \begin{equation} \label{f2} \begin{array}{ccl}
 \lambda_{i,j} &=& -2 e^{2\phi}\phi_{j} \phi_s \bar g^{s p} g_{p i}-e^{2\phi}\phi_{s,j} \bar g^{s p} g_{p i}+e^{2\phi}\phi_s  \bar g^{s q} \bar g_{q \ell,j} \bar g^{\ell p} g_{p i} \\
 &\stackrel{\eqref{LC}}{=}&   -e^{2\phi} \bar g^{s p} g_{p i}(\phi_{s,j} -  \phi_s\phi_j)+e^{2\phi}\phi_s \phi_p \bar g^{s p}  g_{ i j }    \\
 &\stackrel{\eqref{f1}}{=}& -e^{2\phi} \bar g^{s p} g_{p i}
 \left(\tfrac{R}{2(n-1)} g_{s j}-\tfrac{\bar R}{2(n-1)} \bar g_{s j}\right)
 +e^{2\phi}\phi_s \phi_p \bar g^{s p}   g_{ i j }, 
    \end{array}     
 \end{equation} 
 where $\bar g ^{sp}$ is the tensor dual to $\bar g_{i j}$.  
We  combine this with \eqref{a}  and see that 
\begin{equation}\label{vb} 
 \lambda_{i,j}= \mu g_{ij}+ B a_{ij},  
 \end{equation}
 where   $B:=-\tfrac{R}{2(n-1)}$ and  
 $\mu:= \tfrac{\bar R}{2(n-1)}e^{2\phi}+ e^{2\phi}\phi_s \phi_q \bar g^{s q}   $.

Note that $B= -\tfrac{R}{2(n-1)}$ is constant if and only if the scalar curvature $R$ is a constant.

For further use let us also consider the $(1,3)$-tensor 
\begin{equation}\label{X} 
X^i_{ \ jk\ell}:= R^i_{ \ jk\ell} +   \tfrac{ R}{2(n-1)}\left( \delta^i_{\ \ell} g_{jk} - \delta^i_{\ k} g_{j\ell}\right) . \end{equation} 
This tensor clearly satisfies the same algebraic symmetries w.r.t. $g$ as the curvature tensor; by construction of $X$   the contraction 
 $X^s_{\ j s k}= R_{jk} -  \tfrac{R}{2} g_{jk}$ is the stress energy tensor of $g$.  Let us observe  that $X^i_{\ jk\ell}$ satisfies 
 \begin{equation} \label{intX}
 a_{si}X^s_{\ j k\ell}+ a_{sj}X^s_{\ i k\ell}=0:
 \end{equation}
indeed, we substitute $\lambda_{i,j} $ given by \eqref{vb} with   $B= -\tfrac{R}{2(n-1)}$ in \eqref{int1} and  obtain \eqref{X}.

   \subsection{Proof of Theorem \ref{main} under the assumption of Case 1: $R=\const \ne 0$. } \label{case1}

   Without loss of generality we can assume $\frac{R}{2(n-1)}=1,  $ since it always can be achieved by the rescaling of the metric. 
  In this setting 
the  system~\eqref{eq:mg+Ba}  reads
\begin{gather}
  \label{eq:mg-a}
  \left\{
    \begin{array}{ll}
      a_{i j,k}&=\lambda_i g_{j k} + \lambda_j g_{i k}\\
      \lambda_{i,j}&=\mu g_{i j}-a_{i j}\\
      \mu_{,i}&=-2\lambda_i
    \end{array}
  \right.
\end{gather}

By the  
\emph{metric cone} over  $(M,g)$ we understand  the  product manifold $\wt
M=\R_{{}>0}(r)\times M(x)$ equipped by  the  metric $\tg$ such that  in the
coordinates $(r,x)$ its matrix has the form
\begin{gather}
  \label{eq:conemetric}
  \tg(r,x)=\left(
    \begin{array}{cc}
      1& 0 \\
      0& r^2 g(x)
    \end{array}
  \right).
\end{gather}

  Let us recall the following relation between the parallel symmetric $(0,2)$-tensors on the cones and the solutions of \eqref{eq:mg-a}. 
  
\begin{Th}[ Proposition 3.1 and Corollary 3.3  of \cite{Mounoud2010}]
  \label{thm:parallel}
  If a symmetric tensor field $a_{ij}$ on  $(M,g)$  satisfies  \eqref{eq:mg-a}, 
  then the $(0,2)$-tensor  field $A$ on  $(\wt M, \tg)$ defined  in the local coordinates $(r,x)$ 
  by the following (symmetric) matrix:
  \begin{gather}
    \label{eq:A}
    A = \left(
      \begin{array}{c|ccc}
        \mu(x) & -r\lambda_1(x) &\dots &-r\lambda_n(x)\\
        \hline 
        -r\lambda_1(x)& & & \\
        \vdots& & r^2 a(x)& \\
        -r\lambda_n(x)& & & 
      \end{array}
    \right),
  \end{gather}
  is parallel  with respect to the Levi-Civita connection
  of $\tg$.

 Moreover, if a symmetric $(0,2)$-tensor $A_{ij}$ on $\wt M$ is parallel, then in the  cone coordinates it has the
  form~\eqref{eq:A}, where  $(a_{ij}, \lambda_i, \mu)$ 
  satisfy~\eqref{eq:mg-a}.
\end{Th}

\begin{Rem} Since  Proposition 3.1 and Corollary 3.3  of \cite{Mounoud2010}  are written in different mathematical language, let us note that  the proof of  Theorem \ref{thm:parallel}   is actually  an easy exercise. A straightforward way  to do this exercise is to calculate the Levi-Civita connection of the metric $\wt g$ (was done many times before),
  to write down   the condition that a symmetric $(0,2)$-tensor field  on the cone is parallel, and  to compare it with \eqref{eq:mg-a}.
\end{Rem}

    Our next goal is to show (using  the results of \cite{fedorova}) that  the existence of a parallel  symmetric $(0,2)$-tensor field on $\wt M$ 
     that is not proportional to the metric   implies the  existence of a nontrivial parallel 1-form.  We will essentially use that   $n=3 \textrm{ or } 4$;  there are counterexamples to this   claim   in all higher  dimensions (see eg \cite[\S 3.3.3]{fedorova} for a counterexample in dimension $n+1=6$).  
     
     Since $ n=dim(M)\le 4$, the dimension of $\wt M$  is $n+1 \le 5$. Then, the signature of $\wt g$ is as required in the assumptions of \cite[Theorems 5,6]{fedorova}. Then, by \cite[Theorems 5,6]{fedorova}, the dimension of the space of  symmetric $(0,2)$-tensor fields is $\frac{k(k+1)}{2} + \ell$, where $\ell = 1,...,  \lfloor  \frac{n+1-k}{3} \ \rfloor $, where $k$ is the dimension of the space of parallel vector fields and    the brackets ``$\lfloor \ , \ \rfloor$'' mean the
integer part.      Now, for $n=3, 4$ we evidently have  $\lfloor  \frac{n+1-k}{3} \ \rfloor\le 1$. Then, the existence of a 
parallel  symmetric $(0,2)$-tensor field that is not proportional to the metric   implies that $k\ge 1$, i.e., the existence of a nontrivial parallel vector field, which implies the existence of a nontrivial parallel 1-form.

We will call this parallel  1-form by  $V_\alpha$ ($\alpha= 0,...,n$); we work in the cone coordinates $(x_0:= r, x_1,...,x_n)$; we  will denote  the $0$-component  of   $V$ by $v$ so the 1-form $V$ has entries  
$(v, V_1,..., V_n)$. \weg{It is easy to see that $v$ can be viewed as  a function on $M$ and $v_i=(v_1,...,v_n)$ is a 1-form on $M$.}

Since $V_\alpha$ is $\wt g$-parallel, the $(0,2)$-tensor field $V_\alpha V_\beta$  is  also parallel (on $\wt M$);  in the cone  coordinates $x_0:= r,x_1,...,x_n$ 
it is given by the matrix \begin{gather}
    \label{eq:A1}
    A = \left(
      \begin{array}{c|ccc}
        (v)^2  & v V_1 &\dots &vV_n\\
        \hline 
        vV_1& & & \\
        \vdots& &  V_iV_j& \\
        vV_n& & & 
      \end{array}
    \right).
  \end{gather}

Comparing this with  \eqref{eq:A}, we see that by   Theorem \ref{thm:parallel}  $v$ does not depend on $x_0$ (so it is essentially a function on $M$); $V_i$ for $i\ge 1$ have the form $V_i= r v_i$, where $v_i=(v_1,...,v_n)$ is a 1-form on $M$. Moreover, 
the triple $(a_{ij}= v_i v_j , \lambda_i = - v v_i, \mu = (v)^2)$ is a solution of \eqref{eq:mg-a}. 
Note that $(v_1,...,v_n)$ is not zero at a generic point of $M$ since otherwise the `cone' vector field $\tfrac{\partial  }{\partial x_0} $  will be proportional to $V^\alpha$ (with a possible functional coefficient of the proportionality) which is impossible  by \cite[Lemma 4]{fedorova}.
 Note that 
the last  equation of the system \eqref{eq:mg-a} for this solution looks 
$     \left(v^2\right)_{,i}=2 v v_i \textrm{\ implying $v_i=v_{,i}.$ }
  $
  Combining this with $(a_{ij}= v_i v_j , \lambda_i = - v v_i, \mu = v^2)$ and with $\lambda_{i,j}= \mu g_{ij}  - a_{ij}$, we obtain \begin{equation} \label{eq} v_{i,j}= - v g_{ij}.\end{equation} 
  
 Since the matrix of $a_{ij}= v_iv_j$ has rank two, its minimal polynomial has degree two. Then, by Lemma \ref{lem1},  in a neighborhood of a generic point there exists a coordinate system such that the metric $g$ and  the (1,1)-tensor $a^{i}_{\ j}$ 
 are  given by    \eqref{wp}  (cf \cite[Lemma 2.1]{kuhnel}).  Then, the constand $C$ in   \eqref{wp} is $0$  and the 1-form  $v_i$ is given by $v_i=(\sqrt{\pm \sigma}, 0,...,0)$. We see that in this coordinate system $v^i$  and $\lambda^i$ are  proportional to $\tfrac{\partial }{\partial x_1}$ and the coefficient of the proportionality depends on $x_0$.

Next, let us observe that  the vector field $\lambda^i$ is  a projective vector field of $g$, that is, the pullback of $g$ w.r.t. the (local) flow of $\lambda^i$ is geodesically equivalent to $g$. 
  Indeed, it is known (see eg \cite{archive,s1}), that a vector field $v$ is projective if and only if for  the Lie derivative 
  $\ell_{ij}:= \mathcal{L}_v g  $     the $(0,2)$ tensor 
  $a^v:= \ell_{ij}- \tfrac{1}{n+1} \ell^{s}_{\ s} g_{ij}$   satisfies \eqref{basic}.  For $v^i= \lambda^i$, 
  the tensors $\tfrac{1}{2} \ell_{ij}$ and  $\tfrac{1}{2}(\ell_{ij}- \tfrac{1}{n+1} \ell^{s}_{\ s} g_{ij})$ are given by 
  $$
\tfrac{1}{2}  \ell_{ij}= \lambda_{i,j} = \mu  g_{ij} - a_{ij}, \   \tfrac{1}{2}(\ell_{ij}- \tfrac{1}{n+1} \ell^{s}_{\ s} g_{ij})= \mu g_{ij}  -a_{ij} - \tfrac{1}{n+1}  (n\mu - 2 \lambda) g_{ij} = -a_{ij} + \textrm{Const} g_{ij};   
  $$ 
  in the last equality in the  formula above we use that $\mu_i= 2 v v_i= - 2 \lambda_i$, so $\mu= -2 \lambda  + \const$.   Thus, $\ell_{ij}- \tfrac{1}{n+1} \ell^{s}_{\ s} g_{ij}$   satisfies \eqref{basic} and $\lambda^i$   is a projective vector field. 
  
Take a small time $t$ and   denote by $\bar g$  the pullback of $g$ with respect to the time-$t$-flow  of $\lambda_i$. Since  as we explained above in the  coordinate system $(x_1,... , x_n)$  constructed above 
the vector field $\lambda^i$ has the form $(\lambda^1(x_1),0,...,0)$, the pullback of the metric $g$ has also the warped product form $\alpha(x_1) dx_1^2 + \beta(x_1) h(x_2,...,x_n)$, where $h=h_{ij},$ $i,j=2,...,n$,  is a metric on $U\subset  \R^{n-1}(x_2,...,x_n)$. 
Then, the minimal polynomial of  the correspondent  $A$  has degree 2 and by Lemma \ref{lem1} in a certain coordinate system the metrics  $g$ and $\bar g$  have  the form \eqref{wp},  \eqref{barwp}.

Now, since the metric $\bar g$ is isometric to the metric $g$ (since it is its pull-back), every solution  
of equation \eqref{basic}  satisfies \eqref{vn} with $B= -\tfrac{R}{2(n-1)}$. Indeed, both $R$  and $B$ are invariants in the sense they do not  depend  on the coordinate system; moreover, as we explained above, from the result of \cite{KM} it follows that the constant $B$ is the same for all solutions. 
 Then, by Lemma \ref{const},  it satisfies \eqref{f1}, which in view of  $B= \bar B=-\tfrac{R}{2(n-1)}=-\tfrac{\bar R}{2(n-1)}$ reads 
\begin{equation}\label{f20} 
 \phi_{i,j} - \phi_{i} \phi_{j} = \tfrac{R}{2(n-1)}  g_{ij}-\tfrac{R}{2(n-1)} \bar g_{ij} = \tfrac{R}{2(n-1)}  g_{ij}-\tfrac{\bar R}{2(n-1)} \bar g_{ij}  .
 \end{equation} 
Substituting this in   \eqref{ric}, we obtain  \eqref{f11}, which is equivalent to \eqref{stress}. Now, by Lemma \ref{lem1},  geodesically equivalent metrics of the form \eqref{wp}, \eqref{barwp} are affinely equivalent. 
Then,  $R= 0$  which contradicts the assumptions.

 \subsection{Proof of Theorem \ref{main} under the assumption  $dR\ne 0$. } \label{case2}

 We assume as usual that $g$  and $\bar g$ on $M^n$ with $n\ge 3$  are geodesically equivalent and satisfy \eqref{stress}; we show that the assumption 
  that the differential $R_{,i}$   of the scalar curvature $R$ is not zero at a certain point leads to a constradiction.  
  
 As we have shown above, the solution $(a_{ij}, \lambda_i)$ of \eqref{basic} constructed by \eqref{a}, \eqref{lambda} satisfies \eqref{vn}. From the results of  \cite{KM} it follows then that 
 the minimal polynomial of $a^i_{\ j}$ has degree (at most) 
 two.    
 
 More precisely, by \cite[\S 2.3.3]{KM},  under the assumption that a solution $(a_{ij}, \lambda_i)$ of \eqref{basic} satisfies \eqref{vn} (cf \cite[Eq. (38)]{KM}; in our setting  $B$ from \cite[Eq. (38)]{KM} equals 
 $ -\frac{R}{2(n-1)}$), the formula 
  \cite[Eq. (45)]{KM} holds. If we take  a vector field $\xi^i$ such that $\tfrac{1}{(n-1)2} \xi^iR_{, i}=1$ and contract it with  \cite[Eq. (45)]{KM},   we obtain  \eqref{minimal}.

 Then, by Lemma \ref{lem1}, in a certain coordinate system (in a neighborhood of almost every point)  the metrics $g$ and $\bar g$ are as in \eqref{wp},  \eqref{barwp}. By Lemma \ref{lem2}, the metrics are affinely equivalent
 which implies that $R=0$ which contradicts the assumptions.

\subsection{ Proof in the  case $R=0$.} \label{case3}

We assume that $g$ and $\bar g$ on a connected $M^n$ of dimension $n=4$   are geodesically equivalent, are not affinely equivalent, and satisfy \eqref{stress}.
The proof in dimension $n=3$ is similar, is much easier, and will be left to the reader.  We assume $R=0$.
If $\bar R\ne 0$, then we swap $g$ and $\bar g$ and come to the situation considered in \S \ref{case1} and \S \ref{case2}. We  can  therefore assume  without loss of generality that  $\bar R=0$.

First let us show the (local) existence of a novtrivial parallel 1-form  proportinal to $\phi_i$ (the coefficient of proportionality is a function).  In view of $R= \bar R=0$, the equation \eqref{stress} reads $R_{ij}= \bar R_{ij}$. 
Then,   \eqref{ric} implies  $\phi_{i,j}- \phi_i\phi_j=0$. Recall that 
$\phi_i= \phi_{,i}$ for the function $\phi$ given by \eqref{phi}.  Then, for the $1$-form 
$e^{-\phi} \phi_i$ we have 
$$
(e^{-\phi} \phi_i)_{,j}= -e^{-\phi} \phi_i \phi_{j} + e^{-\phi} \phi_{i,j}= 0.
$$

Now let us show that $A=a^i_{\ j}$ has precisely  one nonconstant eigenvalue, moreover, the algebraic multiplicity of this eigenvalue is one (at a generic point).   In order to prove this claim, let us observe that  
 the tensor $X^i_{\ jk\ell}$  given by \eqref{X} coincides  in view of $R=0$ with $R^i_{\ jk\ell}$ so \eqref{intX} implies 
\begin{equation}\label{nullity}
a_{is}R^s_{\ jk\ell}  + a_{js}R^s_{\ ik\ell} = 0.
\end{equation}
Now, the equation \eqref{nullity} can be equivalently reformulated  as follows:
for every $X=X^k$, $Y=Y^\ell$ the endomorphism 
$$
Z:= R(X,Y)= R^i_{ \ j k\ell} X^k Y^\ell$$
commutes with $A$, i.e., $AZ= ZA$.

Since the metrics are not affinely equivalent, at least one of the eigenvalues of $A$ (considered as a function on $M$) is not constant. 
Assume now there exist two nonconstant eigenvalues of $A$, or the algebraic multiplicity of a nonconstant  eigenvalue  is greater than one. 
Recall that by \cite[Proposition 1]{local}, the geometric multiplicity of  every nonconstant  eigenvalues  is one (at the generic point). Now, it is a standard exercise in linear algebra 
 to check that if $g$-selfadjoint $A$ commutes with $g$-scewselfadjoint $Z$, then every  vector from the generalised eigenspace  of $A$ 
 corresponding  to an eigenvalue of geometric multiplicity one  lies in the kernel of $Z$. 
 Thus, if $A$ has   two nonconstant eigenvalues, or if the  algebraic multiplicity of a nonconstant  eigenvalue  is greater than one, 
 then there exist  two  linearly independent vectors $u= u^ i$ and $v= v^i$  such that the restriction of $g$ to $\textrm{span}(v, u)$ is nondegenerate and  both vectors lie in the kernel of each $Z= R(X,Y)$ implying 
  \begin{equation}\label{nullity1}
u^s R^i_{\ sk\ell}  = v^s R^i_{\ sk\ell} = 0. 
\end{equation}

Take the basic such that the  first two vectors of this basis are $v$ and $u$ and the last two vectors are orthogonal to the first one. 
 In this basis, 
the components of curvature tensor with lowed indexes $R_{ijk\ell}$  are zero, if one of the indexes  $i,j,k,\ell$ is 1 or 2. In view of the algebraic  symmetries of the curvature tensor, we obtain that the only nonzero components of $R_{ijk\ell} $ are 
$$
R_{4343 }= R_{3434}= -R_{4334}= -R_{3443}.$$
Calculating the scalar curvature, we obtain 
$$
R= g^{ik}g^{j\ell} R_{ijk\ell}= 2(g^{33 } g^{44} - g^{34} g^{43}) R_{4343}= 2 \det\begin{pmatrix} g^{33} & g^{34} \\ g^{43} & g^{44} \end{pmatrix}R_{4343 } .$$
Since the restriction of $g$ to $ \textrm{span}(v, u)$ is nondegenerate and the first two vectors of our basis are orthogonal to the second two vectors, then 
the determinate in the formula above is not zero so our assumption that the scalar curvature is zero implies that  the curvature tensor  $R^i_{ \ j k\ell} $ vanishes and the metric is flat which contradicts the assumptions.

Thus, only one eigenvalue of  $A$ is not constant, and the geometric multiplicity of this eigenvalue is one. Then, by \eqref{lam}, this eigenvalue is equal to 
$\frac{1}{2} \lambda+ \const$ so $\lambda_i$ is proportional to  this eigenvalue. By the Splitting Lemma  (see \cite[Theorem  3]{splitting}), $\lambda^i$ is an eigenvector 
corresponding to this eigenvalue.  Combining  this with \eqref{lambda}, we see that $\phi_i$ is proportional  to $\lambda_i$.

As we explained in \S \ref{standard},  from  \eqref{vn} and the assumption  $R=0$ it follows the existence of a function $\mu$ such that \eqref{eq:mg+Ba} with $B=0$ holds. 
The third  equation of \eqref{eq:mg+Ba} implies that $\mu= \const$. By \cite{fedorova}, $\mu=0$: indeed,  the second equation of \eqref{eq:mg+Ba} in view $B=0$ reads $\lambda_{i,j}=\mu g_{ij}$.  By scaling the metric we can 
achieve  $\lambda_{i,j}= g_{ij}$.  Now, in   \cite[Lemma 2 and Remark 2]{fedorova} it was shown that an one-form 
$v_i$ such that $v_{i,j}= g_{ij}$ can not be proportional to a parallel one-form at points where the curvature is not zero. Thus, $\mu=0$. 

Since  $\mu=0$, $\lambda_i$ is parallel and therefore is proportional to $e^{-\phi} \phi_i$ with a (nonzero) coefficient. Swamping the metrics $g$ and $\bar g$, we also see that the anolog of the constant $\mu$ (which we, as in Lemma \ref{const}, denote $\bar \mu$, is also  zero). 

Our metrics satisfy the assumptions of Lemma \ref{const} with $B= \bar B= 0$. Then, by Remark \ref{remnew},  in view of $B= \bar \mu=0$, 
we have $g^{pq}\phi_p\phi_q=0$. But $\phi^i$ is an eigenvector of $g$-selfadjoint $A$ such that the corresponding eigenvalue has   algebraic multiplicity $1$, so  $g^{pq}\phi_p\phi_q=0$ implies $\phi_i=0$. Finally, the  metrics are  affinely equivalent.

\section{Counterexample in dimensions $> 4$ 
and   proof of  Corollaries 1,2,3. } 

\subsection{Counterexample}

Consider any metric $h= h_{ij}$ on $U\subseteq \R^{n-1}(x_3,...,x_n)$ of zero scalar curvature. Now, consider a metric 
$g_{ij}$  and a $(1,1)$-tensor $A= a^i_{\ j} $ given by
 \begin{equation} \label{counterexample} 
      g_{ij}= \left(\begin{array}{cc|ccc}
          & x_1  &&&\\ 
         x_1& &&& 
          \\ 
        \hline
        &&&& \\
        &&& (x_2- C)^2h &  \\
        &&&&
      \end{array}\right), \  \ 
     A=  a_{\ j}^i=  \left(\begin{array}{cc|ccc}
          x_2 &x_1 &&&\\     &x_2 &&&\\
        \hline
         && C & & \\
        &&& \ddots &  \\
        &&&& C
      \end{array}\right),
       \end{equation} 
where $C$ is a constant.  By direct calculations one can check that    $g$ is geodesically equivalent to $\bar g$ given by \eqref{LG}. 
Moreover,  the scalar curvatures of $g$ and of $\bar g$ are  equal to zero,  
 and the Riemannian  curvatures of the metrics $g$ and $\bar g$ coincide. Then, the metrics satisfy \eqref{stress}. 
 
 Note that in dimension $n=4$ the metric $h$  is the  flat metric and therefore the metric $g$ is also a flat metric.

\subsection{ Corollaries 1,2,3  follow from \cite{KM, Mounoud2010} }

Corollaries  \ref{cor1}  and   \ref{cor2} easily follow  from the results of our paper combined with that of  \cite{KM}. There, it was proved (under the assumption that the degree of mobility of $g$ is $\ge 3$), that any solution of \eqref{basic} satisfies  \eqref{vn}. Then, it was proved (see \cite[\S\S 2.4, 2.5]{KM}) that the metrics are affinely equivalent provided $g$ is complete and $\bar g$ is light-line complete   (in the indefinite signature) or both metrics are complete (in the definite signature). In the proof the assumption that the degree of mobility is $\ge 3$ was not used  (only \eqref{vn} is necessary for the proof).  

Now, Corollary \ref{cor3} follows from \cite{KM,Mounoud2010}. More precisely,   \eqref{eq:mg+Ba}  implies that the function $\mu$ satisfies the Gallot-Tanno equation 
$$
\lambda_{,ijk} = B (\lambda_i g_{jk} + \lambda_j g_{ik} + 2 \lambda_kg_{ij}),
$$ 
see \cite[Corollary 4]{KM}.
In the case the metrics are not affinely equivalent, $\lambda$ is not constant. Now, by \cite[Theorem 1]{Mounoud2010}, 
the existence of a nonconstant  solution of  this equation implies that the metric has constant  nonzero  sectional curvature. Then,  \eqref{stress} implies that the metrics are proportional.

{\bf Acknowledgements.}
We  thank    Deutsche Forschungsgemeinschaft   and FSU Jena for partial financial support and  H. Bray,  R. Geroch,  G.  Gibbons, D.  Guillini, G. Hall  and  
A.  Wipf  for useful discussions.

\end{document}